\documentclass[a4 paper,11pt,bezier]{article}
\usepackage[top=3cm,bottom=3cm,left=2.5cm,right=2.5cm]{geometry}
\usepackage{amsmath}
\usepackage{amsfonts,amsthm,amssymb}
\usepackage{amsfonts}
\usepackage{graphics,subfigure,caption2}
\usepackage{tikz, color}
\usepackage{graphics,epstopdf}
\usepackage{latexsym,bm}
\usepackage{amsfonts,amsthm,amssymb,mathrsfs,bbding}
\usepackage{CJK}
\usepackage{indentfirst}
\newtheorem{lem}{Lemma}[section]
\newtheorem{thm}[lem]{Theorem}

\theoremstyle{plain}

\begin{document}
\begin{CJK*}{GBK}{song}
\title{Characterizing the forbidden pairs for graphs to be super-edge-connected\footnote{The research is supported by National Natural Science Foundation of China (12261086).}}
\author{Hazhe Ye, Yingzhi Tian\footnote{Corresponding author. E-mail: tianyzhxj@163.com (Y. Tian).}\\
{\small College of Mathematics and System Sciences, Xinjiang
University, Urumqi, Xinjiang 830046, China}}
\date{}
\maketitle
		
\maketitle {\flushleft\bf Abstract:} Let $\mathcal{H}$ be a set of given connected graphs. A graph $G$ is said to be $\mathcal{H}$-free if $G$ contains no $H$ as an induced subgraph for any $H\in \mathcal{H}$. The graph $G$ is super-edge-connected if each minimum edge-cut isolates a vertex in $G$. In this paper, except for some special graphs, we characterize all forbidden subgraph sets $\mathcal{H}$ such that every $\mathcal{H}$-free is super-edge-connected for $|\mathcal{H}|=1$ and $2$.

\maketitle {\flushleft\textit{\bf Keywords}:} Forbidden subgraphs; Super-edge-connectedness; Edge-connectivity
		
\section{Introduction}

In this paper, we consider only finite simple graphs. For notations and graph-theoretical terminology not defined here, we follow \cite{Bondy}.

Let $G=(V,E)$ be a finite simple graph, where $V=V(G)$ is the vertex set and $E=E(G)$ is the edge set.  For a vertex $u\in V(G)$, the $neighborhood$ of $u$ in $G$ is  $N_{G}(u)=\{v\in V(G)|\ v$ is adjacent to $u \}$, and the $degree$ of $u$ in $G$ is $d_G(u)=|N_{G}(u)|$. The $minimum$ $degree$ and $maximum$ $degree$ of $G$ are denoted by $\delta(G)$ and $\triangle(G)$, respectively.  For a vertex set $A\subseteq V(G)$, the $induced$ $subgraph$ of $A$ in $G$, denoted by $G[A]$,  is the graph with vertex set $A$ and two vertices $u$ and $v$ in $A$ are adjacent if and only if they are adjacent in $G$. The $distance$ $d_{G}(u,v)$ between two vertices $u,v\in V(G)$ is the length of the shortest path between  them in the graph $G$.

An edge set $F\subseteq E(G)$ is called an $edge$-$cut$ if  $G-F$ is disconnected. The $edge$-$connectivity$ $\kappa^{\prime}(G)$ of a  graph $G$ is defined as the minimum cardinality of an edge-cut over all edge-cuts of G. The $vertex$-$connectivity$ $\kappa(G)$ of $G$  can be defined similarly. The graph $G$ is connected if $\kappa^{\prime}(G)\geq1$. If $G$ is not connected, them each maximal connected subgraph is call a $component$ of $G$.
It is well known that $\kappa (G)\leq\kappa^{\prime}(G)\leq\delta(G)$. If $\kappa^{\prime}(G)=\delta(G)$, then the graph $G$ is said to be $maximally$ $edge$-$connected$. In 1981, Bauer,  Suffel,   Boesch, and Tindell \cite{Bauer} proposed the concept of super-edge-connectedness. A graph $G$ is said to be $super$-$edge$-$connected$ if each minimum edge-cut isolates a vertex of $G$, that is,  each minimum edge-cut  of $G$ is the  set of edges incident with some vertex in $G$. Clearly, if $G$ is super-edge-connected, then it is also maximally edge-connected; the converse is not true, for example, the cycle $C_n$ ($n\geq4$) is maximally edge-connected but not super-edge-connected. There are lots of sufficient conditions for a graph to be maximally edge-connected, such as minimum degree condition \cite{Chartrand}, Ore condition \cite{Lesniak} and diameter condition \cite{Plesnik}.  Similar sufficient conditions for a digraph to be super-edge-connected can be seen in \cite{Fiol}. For more results, please refer to the survey paper \cite{Hellwig} and the references herein.

Let $\mathcal{H}$ be a set of given connected graphs. A graph $G$ is said to be $\mathcal{H}$-$free$ if there is no induced subgraph in $G$ isomorphic to some graph $H\in \mathcal{H}$. Each graph $H$ in $\mathcal{H}$ is called a $forbidden$ $subgraph$. If $|\mathcal{H}|=2$, then the two graphs in $\mathcal{H}$ are called the $forbidden$ $pair$. For two sets of connected graphs $\mathcal{H}_{1}$ and $\mathcal{H}_{2}$, if for each graph $H_{2}$ in $\mathcal{H}_{2}$, there exists a graph $H_{1}$ in $ \mathcal{H}_{1}$ such that $H_{1}$ is an induced subgraph of $H_{2}$, then denote this relation by  $\mathcal{H}_{1}\preceq \mathcal{H}_{2}$. Obviously,  $\mathcal{H}_{1}\preceq \mathcal{H}_{2}$ implies that a $\mathcal{H}_{1}$-free graph is also $\mathcal{H}_{2}$-free.  Figure 1  gives some classes of forbidden subgraphs used in this paper.

\begin{figure}[htbp]
 % \vspace{-1.0em}
  \centering
  \includegraphics[width=10cm]{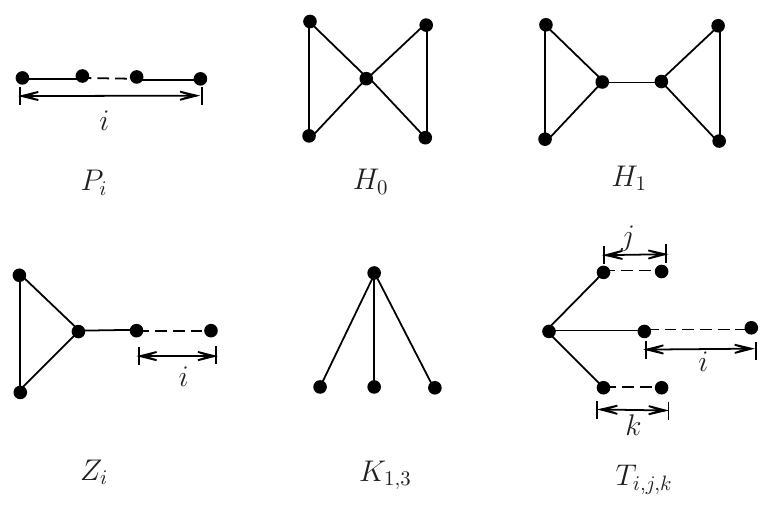}\\
  \caption{Some classes of  forbidden subgraphs.}
  \label{1}
\end{figure}

In \cite{Faudree},  Faudree and Gould determined  the forbidden pairs for hamiltonicity of 2-connected graphs except for a finite number of exceptions.
Wang, Tsuchiya and Xiong \cite{Wang} characterized all the forbidden pairs $R,S$ such that every connected $\{R,S\}$-free graph $G$ has $\kappa(G)=\kappa^{\prime}(G)$. 
Recently,  Du, Huang and Xiong \cite{Du} characterized all the forbidden pairs $R,S$ such that every connected $\{R,S\}$-free graph $G$ is maximally edge-connected.

\begin{thm} (\cite{Wang}) 
 Let $S$ be a connected graph. Then $G$ being a connected $S$-free graph implies $\kappa(G)=\kappa^{\prime}(G)$ if and only if $S$ is an induced subgraph of $P_3$.
\end{thm}

\begin{thm} (\cite{Wang}) 
Let $\mathcal{H}$ be a set of two connected graphs such that each member of $\mathcal{H}$ is not an induced subgraph of $P_3$. Then $G$ being a connected $\mathcal{H}$-free graph implies $\kappa(G)=\kappa^{\prime}(G)$ if and only if $\mathcal{H} \preceq\left\{Z_1, P_5\right\}, \mathcal{H} \preceq\left\{Z_1, K_{1,4}\right\}, \mathcal{H} \preceq\left\{Z_1, T_{1,1,2}\right\}$, $\mathcal{H} \preceq\left\{H_0, P_4\right\}$, or $\mathcal{H} \preceq\left\{ H_0,K_{1,3}\right\}$.
\end{thm}

\begin{thm} (\cite{Du}) 
Let $S$ be a connected graph. Then $G$ being a connected S-free graph implies $\kappa^{\prime}(G)=$ $\delta(G)$ if and only if $S$ is an induced subgraph of $P_4$.
\end{thm}

\begin{thm} (\cite{Du}) 
Let $\mathcal{H}=\{R, S\}$ be a set of two connected graphs such that both $R$ and $S$ are not an induced subgraph of $P_4$. Then $G$ being a connected $\mathcal{H}$-free graph implies $\kappa^{\prime}(G)=\delta(G)$ if and only if $\mathcal{H} \preceq\left\{H_1, P_5\right\}, \mathcal{H} \preceq\left\{Z_2, P_6\right\}$, or $\mathcal{H} \preceq\left\{Z_2, T_{1,1,3}\right\}$.
\end{thm}

Motivated by the results above, we will characterize the forbidden pairs for a graph to be super-edge-connected in this paper. In the next section, the main results will be presented. The proof of the main theorem will be given in the last section.

\section{Main results}
			
\begin{thm}\label{theorem 2.1}
Let $S$ be a connected graph. Then $G$ being a connected S-free graph implies $G$ is super-edge-connected if and only if $S$ is an induced subgraph of $P_3$.
\end{thm}

\renewcommand\proofname{\bf{Proof}}
\begin{proof}

Let $G$ be a connected $P_3$-free graph. Then $G$ is a complete graph, and thus $G$ is super-edge-connected. 

Let $S$ be a connected graph such that every connected $S$-free graph is super-edge-connected. Since $G_{i}$ in Figure 2 is not super-edge-connected for $ i\in \{1, \cdots, 6\}$, then we know $G_{i}$ is not $S$-free and contains $S$ as an induced subgraph. We observe that the common induced subgraph of the graphs in ${G}_3$ and ${G}_6$ is a path and the longest induced path of the graphs in ${G}_1$ is $P_3$, then $S$ must be an induced subgraph of $P_3$. The proof is complete.
\end{proof}

The $Cartesian$ $product$ of two graphs $G_1$ and $G_2$, denoted by $G_1 \square G_2$, is defined on the vertex sets $V(G_1) \times V(G_2)$, and $(x_1, y_1)(x_2, y_2)$ is an edge in $G_1 \square G_2$ if and only if one of the following is true: ($i$) $x_1=x_2$ and $y_1 y_2 \in E(G_2)$;
($ii$) $y_1=y_2$ and $x_1 x_2 \in E(G_1)$.

In the following, we try to characterize the forbidden pairs for a graph to be super-edge-connected. By Theorem 2.1,  a connected $P_3$-free graph is super-edge-connected. Thus for any connected graph $S$, we obtain that if $G$ is $\left\{P_{3} , S\right\}$-free, then $G$ is super-edge-connected. So we assume that both $R$ and $S$ are not an induced subgraph of $P_3$ in the following main theorem in this paper.

\begin{thm}\label{theorem 2.2}
Let $\mathcal{H}=\{R, S\}$ be a set of two connected graphs such that both $R$ and $S$ are not an induced subgraph of $P_3$. Then $G$ being a connected $\mathcal{H}$-free graph implies $G$ is super-edge-connected  if and only if (i) $G\ncong C_{4}$ and $\mathcal{H} \preceq\left\{H_0, P_4\right\}$,  or (ii) $G\ncong P_{2}\Box P_{3}$, $G\ncong P_{n}, C_{n}\ (n\geq4)$ and $\mathcal{H} \preceq\left\{Z_1, T_{1,1,2}\right\}$.
\end{thm}

\section{The Proof of Theorem 2.2}
In Figure 2, we construct some classes of non super-edge-connected   graphs $G_{i}$ (where $i\in \{1, \cdots, 6\}$) on $n$ vertices, which will be used in the proof of the main result. To distinguish these graphs, we assume $n\geq8$.  $C_n$ denotes the cycle on $n$ vertices and $K_n$ denotes the complete graph on $n$ vertices.

\begin{figure}[htbp]
  %\vspace{-1.0em}
  \centering
  \includegraphics[width=10cm]{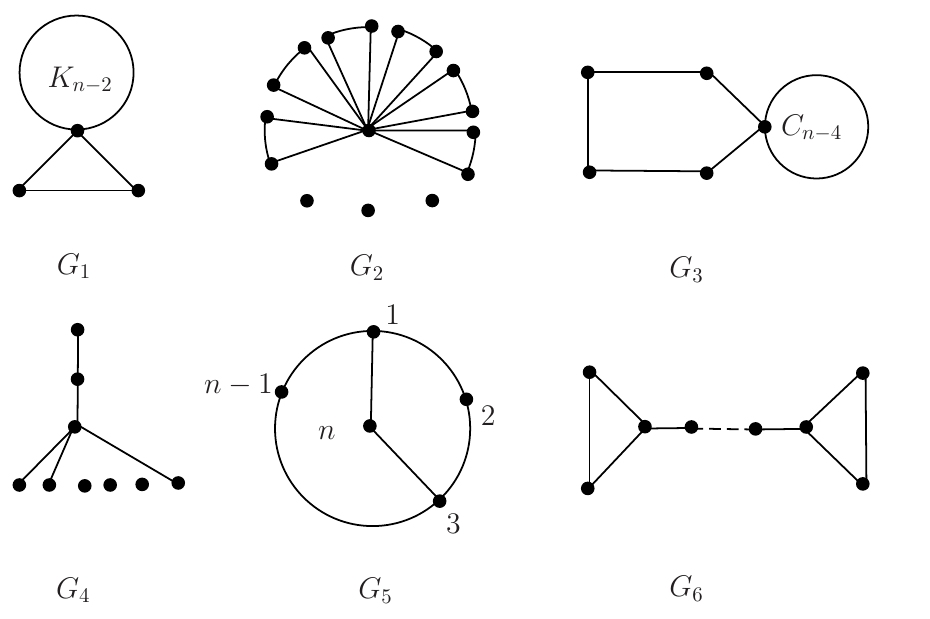}\\
  \caption{Some classes of non super-edge-connected graphs.}
  \label{2}
\end{figure}

\renewcommand\proofname{\bf{The Proof of Theorem 2.2.}}

\begin{proof}
{\bf{We first prove the necessity.}} Assume $\mathcal{H}=\{R,S\}$ is a set of connected graphs such that both $R$ and $S$ are not an induced subgraph of $P_{3}$, and every connected $\mathcal{H}$-free graph is super-edge-connected. Let $G$ be a connected $\mathcal{H}$-free graph. Then $G\ncong P_2\Box P_3$ and $G\ncong P_n,C_n\ (n\geq4)$. Since $G_{i}$ in Figure 2 is not super-edge-connected, we obtain that $G_{i}$ is not $\mathcal{H}$-free  and  contains at least one of $R$ and $S$ as an induced subgraph for $i\in \{1, \cdots, 6\}$.

\noindent\textbf{Claim 1.} Either $R$ or $S$ is an induced subgraph of $H_{0}$.

By contradiction, assume that neither $R$ nor $S$ is an induced subgraph of $H_{0}$. Without loss of generality, assume that $R$ is an induced subgraph of $G_{1}$. Since $R$ is not an induced subgraph of $H_{0}$. We observe that $R$ must contain a $K_{t}$ as an induced subgraph with $t\geq4$ . For $i\in\{2,3,4,5,6\}$, $G_{i}$ is $R$-free. Then they contain $S$ as an induced subgraph. We observe that the common induced subgraph of the graphs in ${G}_3$ and ${G}_6$ is a path and the longest induced path of $G_{2}$ is $P_{3}$, then $S$ should be an induced subgraph of $P_{3}$, a contradiction.  Claim 1 is thus proved.

By Claim 1, assume, without loss of generality, that $R$ is an induced subgraph of $H_{0}$. We distinguish two cases as follows.

\noindent\textbf{Case 1.} $R$ is $H_{0}$.

For $i\in\{3,4,5,6\}$, $G_{i}$ is $R$-free. Then they contain $S$ as an induced subgraph. We observe that the common induced subgraph of the graphs in ${G}_3$ and ${G}_6$ is a path and the longest induced path of $G_{4}$ is $P_{4}$, then $S$ should be an induced subgraph of $P_{4}$. So $\mathcal{H}=\{R,S\}\preceq\left\{H_{0},P_{4}\right\}$ .

\noindent\textbf{Case 2.} $R$ is an induced subgraph of $Z_{1}$.

Since $R$ is not an induced subgraph of $P_3$, we get that $R$  must contain a triangle. For $i\in\{3,4,5\}$, $G_{i}$ is $R$-free. Then they contain $S$ as an induced subgraph. We observe that the maximal common induced subgraph of the graphs in ${G}_4$ and ${G}_5$ is a $T_{1,1,2}$, then $S$ should be an induced subgraph of $T_{1,1,2}$. So $\mathcal{H}=\{R,S\}\preceq\left\{Z_{1},T_{1,1,2}\right\}$.

{\bf{Now we are going to prove the sufficiency.}}
By contradiction, assume $G$ is a connected $\mathcal{H}$-free graph, but $G$ is not super-edge-connected, where ($i$) $G\ncong C_{4}$ and $\mathcal{H} \preceq\left\{H_0, P_4\right\}$,  or ($ii$) $G\ncong P_{2}\Box P_{3}$, $G\ncong P_{n}, C_{n}\ (n\geq4)$ and $\mathcal{H} \preceq\left\{Z_1, T_{1,1,2}\right\}$.  Since $G$ is not super-edge-connected, there is a minimum edge-cut $F$ such that $G-F$ has no isolated vertices. Clearly, $G-F$ has only two components, say $X$ and $Y$.
Let $X_1=V\left(X\right) \cap V(F)$ and $Y_1=V\left(Y\right) \cap V(F)$.  Denote $X_2=V(X)\setminus X_1$ and $Y_2=V(Y)\setminus Y_1$.
Note that $|V(X)|, |V(Y)|\geq2$ and $|X_1|, |Y_1|\leq |F|=\kappa^{\prime}(G)\leq\delta(G)$.

\noindent\textbf{Claim 2.} If $X_2=\emptyset$, then $X$ is a complete graph with order $\delta(G)$, each vertex in $X$ is incident with exactly one edge in $F$, and $|X|=|X_1|=|F|=\kappa^{\prime}(G)=\delta(G)$; similarly, if $Y_2=\emptyset$, then $Y$ is a complete graph with order $\delta(G)$, each vertex in $Y$ is incident with exactly one edge in $F$, and $|Y|=|Y_1|=|F|=\kappa^{\prime}(G)=\delta(G)$.

By $\kappa^{\prime}(G)\leq\delta(G)$, $|X_1|\leq|V(X)|$ and $|X_1|\leq |F|=\kappa^{\prime}(G)$, we have the following inequalities.

$$
\begin{aligned}
\left|E\left(X\right)\right| & =\frac{1}{2}\left(\sum_{x \in V\left(X\right)} d_G(x)-\kappa^{\prime}(G)\right) \\
& \geq \frac{1}{2}\left(\delta(G)\left|V\left(X\right)\right|-\kappa^{\prime}(G)\right) \\
& \geq \frac{1}{2}\left(\delta(G) |X_1|-\kappa^{\prime}(G)\right) \\
& \geq\frac{1}{2} \kappa^{\prime}(G)\left(|X_1|-1\right) \\
& \geq \frac{1}{2} |X_1|\left(|X_1|-1\right).
\end{aligned}
$$

If $X_2=\emptyset$, then all the inequalities above will be equalities. Thus we obtain that  $X$ is a complete graph with order $\delta(G)$, each vertex in $X$ is incident with exactly one edge in $F$, and $|X|=|X_1|=|F|=\kappa^{\prime}(G)=\delta(G)$.

By Claim 2, we consider two cases in the following.

\noindent\textbf{Case 1.} At least one of $X_2$ and $Y_2$ is an empty set.

Assume, without loss of generality, that $X_2=\emptyset$.  Then, by Claim 2,   $X$ is a complete graph with order $\delta(G)$, each vertex in $X$ is incident with exactly one edge in $F$, and $|X|=|X_1|=|F|=\kappa^{\prime}(G)=\delta(G)$. In addition, by $|V(X)|\geq2$, we have $\delta(G)\geq2$.

\noindent\textbf{Subcase 1.1.} $G\ncong C_{4}$ and $\mathcal{H} \preceq\left\{H_0, P_4\right\}$.

Suppose $Y_2=\emptyset$. Then $G$ is the union of two complete graphs on $\delta(G)$ vertices, together with the perfect matching $F$, that is, $G$ is isomorphic to the Cartesian product graph  $K_{\delta(G)}\Box K_2$. If $\delta(G)\geq3$, then there exists an induced $P_{4}$, a contradiction. Otherwise, if $\delta(G)=2$, then $G\cong C_4$, also a contradiction. 

Suppose $Y_2\neq\emptyset$. Then by $Y$ is connected, there is a path $x_{1}x_{2}y_{1}y_{2}$, where $x_{1},$ $x_{2}\in V\left(X\right)$, $y_{1}\in Y_1$ and $y_{2}\in Y_2$. 
If $x_1y_1\notin E(G)$, then $x_{1}x_{2}y_{1}y_{2}$ is an induced path of $G$, a contradiction. Thus we assume $x_1y_1\in E(G)$. 
Since $d_{G}\left(y_{2}\right)\geq\delta(G)\geq2$,  $y_2$ has  a neighbor $y_3$ different from $y_1$. If $y_1y_3\notin E(G)$, then $G\left[\{x_1,y_1,y_2,y_3\}\right]\cong P_4$,  a contradiction. Otherwise, if $y_1y_3\in E(G)$, then $G\left[\{x_1,x_2,y_1,y_2,y_3\}\right]\cong H_0$, also a contradiction.

\noindent\textbf{Subcase 1.2.} $G\ncong P_{2}\Box P_{3}$, $G\ncong P_{n}, C_{n}\ (n\geq4)$ and $\mathcal{H} \preceq\left\{Z_1, T_{1,1,2}\right\}$.

Clearly, $|X_1|=|F|\geq|Y_1|$. If $|X_1|>|Y_1|$, then there exists a vertex in $Y_1$, say $y_1$, has at least two neighbors in $X_1$, say $x_1$ and $x_2$.
By $Y$ is connected, $y_1$ has at least one neighbor, say $y_2$, in $V(Y)$.  Since $X$ is a complete graph and each vertex in $X$ is incident with exactly one edge in $F$, we 
know that  $G\left[\{x_1,x_2,y_1,y_2\}\right]\cong Z_1$, a contradiction. Thus, we assume $|X_1|=|Y_1|$ in the following.

Since $|X_1|=|Y_1|=|F|$, we know that each vertex in $X_1$ has only one neighbor in $Y_1$ and each vertex in $Y_1$ has only one neighbor in $X_1$. 
If $|X_1|=|Y_1|\geq3$, then the induced subgraph by any three vertex in $X_1$, say $x_1,x_2,x_3$, together with the neighbor of $x_1$ in $Y_1$ is isomorphic to $Z_{1}$, a contradiction. It is remaining to consider $|X_1|=|Y_1|=2$.  Assume $X_1=\{x_1,x_2\}$, $Y_1=\{y_1,y_2\}$ and $F=\{x_1y_1,x_2y_2\}$. 

Suppose $y_1y_2\in E(G)$. Since $G$ is not a cycle, $y_1$ or $y_2$, say $y_2$ has a neighbor $y_3$ in $Y_2$. If $y_1y_3\in E(G)$, then $G[\{x_1,y_1,y_2,y_3\}]\cong Z_1$, a contradiction. So $y_1y_3\notin E(G)$. Since $\delta(G)\geq2$, $y_3$ has another neighbor $y_4$ different from $y_2$. If $y_2y_4\in E(G)$, then $G[\{x_2,y_2,y_3,y_4\}]\cong Z_1$, a contradiction. So $y_2y_4\notin E(G)$. If $y_1y_4\notin E(G)$, then $G[\{x_2,y_1,y_2,y_3,y_4\}]\cong T_{1,1,2}$, a contradiction. So $y_1y_4\in E(G)$.  
If $V(G)=\{x_1,x_2,y_1,y_2,y_3,y_4\}$, then $G\cong P_{2}\Box P_{3}$, a contradiction. Thus, by $G$ is connected, there is a vertex $y_5$ adjacent to $y_1$,  $y_2$, $y_3$, or $y_4$. We only consider the case $y_1y_5\in E(G)$ here, others can be analysed similarly. If $y_4y_5\in E(G)$, then $G[\{x_1,y_1,y_4,y_5\}]\cong Z_1$, a contradiction. 
Otherwise, if  $y_4y_5\notin E(G)$, then $G[\{x_1,x_2,y_1,y_4,y_5\}]\cong T_{1,1,2}$, also a contradiction. 

Suppose $y_1y_2\notin E(G)$.  Since $Y$ is connected, there is a path $z_1z_2\cdots z_t$ connecting $y_1$ and $y_2$ in $Y$, where $z_1=y_1$, $z_t=y_2$ and $t\geq3$.
Furthermore, we choose this path  $z_1z_2\cdots z_t$ to be a shortest path connecting $y_1$ and $y_2$ in $Y$. Since $G$ is not  a cycle, there is a vertex $w\in Y\setminus\{z_1,z_2,\cdots,z_t\}$ such that $w$ is adjacent to some vertex in $\{z_1,z_2,\cdots,z_t\}$.  Assume $wz_i\in E(G)$ and $wz_j\notin E(G)$ for $j<i$. Let $z_0=x_1$ and $z_{-1}=x_2$ . If $wz_{i+1}\in E(G)$, then $G[\{w,z_{i-1},z_{i},z_{i+1}\}]\cong Z_1$, a contradiction.  Otherwise, if  $wz_{i+1}\notin E(G)$, then $G[\{w,z_{i-2},z_{i-1},z_{i},z_{i+1}\}]\cong T_{1,1,2}$, also a contradiction.

\noindent\textbf{Case 2.} Neither $X_2$ nor $Y_2$ is an empty set.

Since the distance between a vertex in $X_2$ and a vertex in $Y_2$ is at least three, $G$ must contains $P_4$ as an induced path. So the case $G\ncong C_{4}$ and $\mathcal{H} \preceq\left\{H_0, P_4\right\}$ can not happen. Thus we only consider the case that $G\ncong P_{2}\Box P_{3}$, $G\ncong P_{n}, C_{n}\ (n\geq4)$ and $\mathcal{H} \preceq\left\{Z_1, T_{1,1,2}\right\}$ in the following.  We distinguish  two subcases as follows.

\noindent\textbf{Subcase 2.1.} $G$ contains a path $x_{2}x_{1}y_{1}y_{2}$, where $x_i\in X_i$ and $y_i\in Y_i$ for $i=1,2$.

\noindent\textbf{Subcase 2.1.1.}  $\left|N_{G}\left(x_{2}\right)\cap X_2\right|\geq2$ or $\left|N_{G}\left(y_{2}\right)\cap Y_2\right|\geq2$.

Assume $\left|N_{G}\left(x_{2}\right)\cap X_2\right|\geq2$. Then there exist two vertices $x_{2}^{\prime}$, $x_{2}^{\prime\prime}\in X_2$ such that $x_{2}x_{2}^{\prime}$, $x_{2}x_{2}^{\prime\prime}\in E(G)$. If $x_{2}^{\prime}x_{2}^{\prime\prime}$, $x_{2}^{\prime}x_{1}$, $x_{2}^{\prime\prime}x_{1}\notin E(G)$, then $G\left[\{x_1,x_2,x_2^{\prime},x_2^{\prime\prime},y_1\}\right]\cong T_{1,1,2}$, a contradiction. If at least one of these three edges exists, then $G$ contains an induced subgraph isomorphic to $Z_1$, also a contradiction. For example, if $x_{2}^{\prime}x_{1}\in E(G)$, then $G\left[\{x_1,x_2,x_2^{\prime},y_1\}\right]\cong Z_1$.
The case $\left|N_{G}\left(y_{2}\right)\cap Y_2\right|\geq2$ can be proved similarly.

\noindent\textbf{Subcase 2.1.2.} $\left|N_{G}\left(x_{2}\right)\cap X_2\right|=1$ and  $\left|N_{G}\left(y_{2}\right)\cap Y_2\right|=1$.

Let $x_3\in N_{G}\left(x_{2}\right)\cap X_2$ and $y_3\in N_{G}\left(y_{2}\right)\cap Y_2$. If $x_1x_3$ or $y_1y_3\in E(G)$, then $G$ contains an induced $Z_{1}$. For example, if $x_1x_3\in E(G)$, then $G\left[\{x_1,x_2,x_3,y_1\}\right]\cong Z_1$. Hence, $x_1x_3\notin E(G)$ and $y_1y_3\notin E(G)$. 

Suppose $d_{G}\left(x_{1}\right)\geq3$.  Let $w\in N_{G}\left(x_{1}\right)\setminus \{x_2,y_1\}$. If $w\in X_2$, then $G\left[\{x_1,x_2,y_1,w\}\right]\cong Z_1$ when $wx_2\in E(G)$, and $G\left[\{x_1,x_2,y_1,y_2,w\}\right]\cong T_{1,1,2}$ when $wx_2\notin E(G)$,  contradicts to the assumption.  By a similar argument,  we can obtain contradictions for $w\in X_1$ and $w\in Y_1$. Thus $d_{G}\left(x_{1}\right)=2$. Similarly, $d_{G}\left(y_{1}\right)=2$. Since $G\ncong P_{n}, C_{n}\ (n\geq4)$, then there  exists a vertex in $G$ such that its degree is at least three.  Assume, without loss of generality, $V(X)$ has a vertex with degree at least three. Choose such a vertex $z_1$ such that its distance to $x_1$ is minimum in the component $X$. Let $z_1z_2\cdots z_r(z_r=x_1)$ be the shortest path between $z_1$ and $x_1$. Then $d_G(z_i)=2$ for $2\leq i\leq r$. Let $z_{1}^{\prime}$, $z_{1}^{\prime\prime}\in N_{G}\left(z_{1}\right)\setminus \{z_2\}$.  Then we can find an induced $Z_1$ or $T_{1,1,2}$ according to $z_{1}^{\prime}z_{1}^{\prime\prime}$ is an edge in $G$ or not.

\noindent\textbf{Subcase 2.1.3.} $\left|N_{G}\left(x_{2}\right)\cap X_2\right|=0$ and  $\left|N_{G}\left(y_{2}\right)\cap Y_2\right|\leq1$, or $\left|N_{G}\left(x_{2}\right)\cap X_2\right|\leq1$ and  $\left|N_{G}\left(y_{2}\right)\cap Y_2\right|=0$.

Without loss of generality, assume that$\left|N_{G}\left(x_{2}\right)\cap X_2\right|=0$ and  $\left|N_{G}\left(y_{2}\right)\cap Y_2\right|\leq1$. Since $N_{G}\left(x_{2}\right)\subseteq X_1$ and $|X_1|\leq|F|=\kappa^{\prime}(G)\leq\delta(G)$, we have $d_{G}\left(x_{2}\right)=|X_1|=|F|=\kappa^{\prime}(G)=\delta(G)$. Note that $|Y_1|\leq|F|=|X_1|$.

Suppose $|X_1|>|Y_1|$. Then $\left|N_{G}\left(y_{2}\right)\cap Y_2\right|\geq1$.  Let $y_3\in N_{G}\left(y_{2}\right)\cap Y_2$. If $y_1y_3\in E(G)$, then $G\left[\{x_1,y_1,y_2,y_3\}\right]\cong Z_1$, a contradiction. So assume $y_1y_3\notin E(G)$. If  $\left|N_{G}\left(y_{1}\right)\cap X_1\right|\geq2$,  then we can find an induced $Z_1$ or $T_{1,1,2}$ according to any two vertices in $N_{G}\left(y_{1}\right)\cap X_1$ are adjacent or not.  So assume $y_1$ has exactly one neighbor $x_1$ in $X_1$.
By $|X_1|>|Y_1|\geq1$, we get $|X_1|\geq2$. Since $y_1$ has exactly one neighbor $x_1$ in $X_1$, we have $|Y_1|\geq2$, thus $|X_1|\geq3$. If there is at least one edge in $G[X_1]$, then we can find an induced $Z_1$, a contradiction. Otherwise, we can find an induced $T_{1,1,2}$ rooted at $x_2$, also a contradiction.

Suppose $|X_1|=|Y_1|$. Then by $|X_1|=|F|=\kappa^{\prime}(G)=\delta(G)$,  we know that each vertex in $X_1$ has only one neighbor in $Y_1$ and each vertex in $Y_1$ has only one neighbor in $X_1$. That  is,  $F$ is perfect matching in $G[X_1\cup Y_1]$ connecting $X_1$ and $Y_1$.

Suppose $|X_1|=|Y_1|\geq3$. Then according to there is an edge in $G[X_1]$ or not, we can find an induced $Z_1$ or $T_{1,1,2}$, a contradiction. Suppose $|X_1|=|Y_1|=2$.
By a similar argument as Subcase 1.2, we can obtain a contradiction. Suppose $|X_1|=|Y_1|=1$. Since $G\ncong P_{n}, C_{n}\ (n\geq4)$, then there  exists a vertex in $G$ such that its degree is at least three. By a similar argument as Subcase 2.1.2, we can also obtain a contradiction.

\noindent\textbf{Subcase 2.2.} $G$ contains no path $x_{2}x_{1}y_{1}y_{2}$ satisfying $x_i\in X_i$ and $y_i\in Y_i$ for $i=1,2$.

Let $X_1^1=\left\{x \in X_1: N_G(x) \cap X_2 \neq \emptyset\right\}$ and $Y_1^1=\left\{y \in Y_1: N_G(y) \cap Y_2 \neq \emptyset\right\}$. Denote  $X_1^2=X_1-X_1^1$ and $Y_1^2=Y_1-Y_1^1$. Since $G$ contains no path $x_{2}x_{1}y_{1}y_{2}$ satisfying $x_i\in X_i$ and $y_i\in Y_i$ for $i=1,2$,  we obtain that  $X_1^2\neq\emptyset$, $Y_1^2\neq\emptyset$ and there is no edge connecting $X_1^1$ and  $Y_1^1$. Since $X$ is connected, there is an edge connecting $X_1^1$ and $X_1^2$, say $x_1^1x_1^2\in E(G)$, where $x_1^1\in X_1^1$ and $x_1^2\in X_1^2$. Let $x_2$ be a neighbor of $x_1^1$ in $X_2$ and $y_1$ be a neighbor of $x_1^1$ in $Y_1$. Since $|X_1^1|<\delta(G)$, $x_2$ has a neighbor, say $x_3$, in $X_2$. If $x_1^1x_3\in E(G)$, then $G\left[\{x_1^1,x_2,x_3,y_1\}\right]\cong Z_1$, a contradiction. So assume $x_1^1x_3\notin E(G)$. If $x_1^2y_1\in E(G)$, then $G\left[\{x_1^1,x_1^2,x_2,y_1\}\right]\cong Z_1$, a contradiction. So assume $x_1^2y_1\notin E(G)$. Thus $G\left[\{x_1^1,x_1^2,x_2,x_3,y_1\}\right]\cong T_{1,1,2}$, also a contradiction.

Since all cases lead to contradictions, the proof is thus complete.
\end{proof}

\end{CJK*}
	
\end{document}